\documentclass[a4paper,12pt]{amsart}
\usepackage{amssymb}
\usepackage{ifthen}
\usepackage{graphicx}
\usepackage{color,xcolor}
\nonstopmode
\numberwithin{equation}{section}
\newtheorem{thm}{Theorem}
\newtheorem{cor}{Corollary}
\newtheorem{lem}{Lemma}
\newtheorem{prop}{Proposition}

\newtheorem{conj}{Conjecture}
\newtheorem{prob}{Problem}

\theoremstyle{definition}
\newtheorem{defn}{Definition}
\newtheorem{ca}{Case}

\newtheorem{rem}{Remark}

\newenvironment{pf}[1][]{%
 \vskip 1mm
 \noindent
 \ifthenelse{\equal{#1}{}}%
  {{\slshape Proof. }}%
  {{\slshape #1.} }%
 }%
{\qed\medskip}

\newcounter{alphabet}

\newenvironment{Thm}[1][]{\refstepcounter{alphabet}%
\bigskip%
\noindent%
{\bf Theorem \Alph{alphabet}}%
\ifthenelse{\equal{#1}{}}{}{ (#1)}%
{\bf .} \itshape}{\vskip 8pt}

\newcounter{alphabet2}

\newcommand{\IN}{{\mathbb N}}

\newcommand{\ID}{{\mathbb D}}

\def\be{\begin{equation}}
\def\ee{\end{equation}}

\newcommand{\ben}{\begin{enumerate}}
\newcommand{\een}{\end{enumerate}}

\newcommand{\blem}{\begin{lem}}
\newcommand{\elem}{\end{lem}}
\newcommand{\bthm}{\begin{thm}}
\newcommand{\ethm}{\end{thm}}
\newcommand{\bcor}{\begin{cor}}
\newcommand{\ecor}{\end{cor}}
\newcommand{\beg}{\begin{exam}}
\newcommand{\eeg}{\end{exam}}
\newcommand{\begs}{\begin{examples}}
\newcommand{\eegs}{\end{examples}}
\newcommand{\bdefe}{\begin{defn}}
\newcommand{\edefe}{\end{defn}}
\newcommand{\bprob}{\begin{prob}}
\newcommand{\eprob}{\end{prob}}
\newcommand{\bques}{\begin{ques}}
\newcommand{\eques}{\end{ques}}
\newcommand{\bei}{\begin{itemize}}
\newcommand{\eei}{\end{itemize}}
\newcommand{\bcon}{\begin{conj}}
\newcommand{\econ}{\end{conj}}
\newcommand{\bop}{\begin{op}}
\newcommand{\eop}{\end{op}}

\newcommand{\bas}{\begin{assertion}}
\newcommand{\eas}{\end{assertion}}

\newcommand{\bfa}{\begin{fact}}
\newcommand{\efa}{\end{fact}}

\newcommand{\bca}{\begin{ca}}
\newcommand{\eca}{\end{ca}}

\newcommand{\bst}{\begin{step}}
\newcommand{\est}{\end{step}}

\newcommand{\bsca}{\begin{sca}}
\newcommand{\esca}{\end{sca}}

\newcommand{\bcl}{\begin{cl}}
\newcommand{\ecl}{\end{cl}}

\newcommand{\bmlem}{\begin{mlem}}
\newcommand{\emlem}{\end{mlem}}

\newcommand{\bscl}{\begin{scl}}
\newcommand{\escl}{\end{scl}}

\newcommand{\bcons}{\begin{conjs}}
\newcommand{\econs}{\end{conjs}}

\newcommand{\bprop}{\begin{prop}}
\newcommand{\eprop}{\end{prop}}

\newcommand{\br}{\begin{rem}}
\newcommand{\er}{\end{rem}}
\newcommand{\brs}{\begin{rems}}
\newcommand{\ers}{\end{rems}}
\newcommand{\bo}{\begin{obser}}
\newcommand{\eo}{\end{obser}}
\newcommand{\bos}{\begin{obsers}}
\newcommand{\eos}{\end{obsers}}
\newcommand{\bpf}{\begin{pf}}
\newcommand{\epf}{\end{pf}}
\newcommand{\ba}{\begin{array}}
\newcommand{\ea}{\end{array}}
\newcommand{\beq}{\begin{eqnarray}}
\newcommand{\beqq}{\begin{eqnarray*}}
\newcommand{\eeq}{\end{eqnarray}}
\newcommand{\eeqq}{\end{eqnarray*}}

\newcommand{\ra}{\to}

\newcommand{\ds}{\displaystyle}

\newcounter{minutes}
\divide\time by 60
\newcounter{hours}
\multiply\time by 60 \addtocounter{minutes}{-\time}

\begin{document}

\bibliographystyle{amsplain}
\title [Bohr-type inequalities for harmonic mappings with a multiple zero]
{Bohr-type inequalities for harmonic mappings with a multiple zero at the origin}

\def\thefootnote{}
\footnotetext{ \texttt{\tiny File:~\jobname .tex,
          printed: \number\day-\number\month-\number\year,
          \thehours.\ifnum\theminutes<10{0}\fi\theminutes}
} \makeatletter\def\thefootnote{\@arabic\c@footnote}\makeatother

\author{Yong Huang}
 \address{Y. Huang, School of Mathematical Sciences, South China Normal University, Guangzhou, Guangdong 510631, China.}
 \email{hyong95@163.com}

\author{Ming-Sheng Liu${}^{~\mathbf{*}}$}
 \address{M-S Liu, School of Mathematical Sciences, South China Normal University, Guangzhou, Guangdong 510631, China.} \email{liumsh65@163.com}

\author{Saminathan Ponnusamy}
\address{S. Ponnusamy, Department of Mathematics,
Indian Institute of Technology Madras, Chennai-600 036, India. }
\email{samy@iitm.ac.in}




\begin{abstract}
In this paper, we first determine Bohr's inequality for the class of harmonic mappings $f=h+\overline{g}$ in the
unit disk $\ID$, where either both $h(z)=\sum_{n=0}^{\infty}a_{pn+m}z^{pn+m}$ and $g(z)=\sum_{n=0}^{\infty}b_{pn+m}z^{pn+m}$
are analytic and bounded in $\ID$, or satisfies the condition
$|g'(z)|\leq d|h'(z)|$ in $\ID\backslash \{0\}$ for some $d\in [0,1]$ and $h$ is bounded. In particular,
we obtain Bohr's inequality for the class of harmonic $p$-symmetric mappings.
Also, we investigate the Bohr-type inequalities of harmonic mappings with a multiple zero at the origin and that most of results are
proved to be sharp.
\end{abstract}

\maketitle
\pagestyle{myheadings}
\markboth{Y. Huang, M-S Liu and S. Ponnusamy}{Bohr-type inequalities for harmonic mappings with a multiple zero}

\section{Preliminaries and some basic questions}\label{HLP-sec1}

The classical theorem of Bohr \cite{B1914}, examined a century ago,  generates intensive research activity--what is called
Bohr's phenomena.
Determination of the Bohr radius for analytic functions in a domain \cite{FR}, as well as for analytic functions
from $\ID$ into particular domains, such as the punctured unit disk, the exterior of the closed unit disk,
and concave wedge-domains, has been discussed in the literature \cite{Abu,Abu2,Abu4,Abu3}. See also
the recent survey articles \cite{AAP2016,IKKP2018,KKP2018} and \cite[Chapter 8]{GarMasRoss-2018}.
The interest in the Bohr phenomena was revived in the nineties due to the extensions to
holomorphic functions of several complex variables and to more abstract settings.
For example in 1997, Boas and Khavinson \cite{BK1997} found bounds for Bohr's radius in any
complete Reinhard domains and showed that the Bohr radius decreases to zero as the dimension of the
domain increases. This paper stimulated interests on Bohr type questions in different settings.
For example, Aizenberg \cite{A2000,A2005}, Aizenberg et al. \cite{AAD}, Defant and Frerick \cite{DF}, and
Djakov and Ramanujan \cite{DjaRaman-2000} have established further results on Bohr's phenomena for multidimensional
power series. Several other aspects and generalizations of Bohr's inequality may be obtained from the literature.
For instance, Defant \cite{DFOOS} improved a version of the Bohnenblust-Hille inequality and Paulsen \cite{PS2004} proved a
uniform algebra analogue of the classical inequality of Bohr concerning Fourier coefficients of bounded holomorphic functions in 2004.
In \cite{PPS2002,PS2006}, the authors demonstrated the classical Bohr inequality using different methods of operators.
Abu Muhanna \cite{Abu} and, Kayumov and Ponnusamy \cite{KP2017} investigated Bohr's inequality for the class of analytic functions that are
subordinate to univalent functions and odd univalent functions, respectively.
On the other hand, Ali et al. \cite{AliBarSoly} discussed Bohr's phenomenon for the classes of even and odd analytic functions
and also for alternating series.  In \cite{AlKayPON-19,KSS2017,LPW2020,LP2019}, the authors considered the Bohr radius for the family
$K$-quasiconformal  sense-preserving harmonic mappings and the class of all sense-preserving harmonic mappings, separately.
Recently, the articles \cite{LSX2018,PVW2019,PW2019} presented a refined version of Bohr's inequality along with few other related improved versions
of previously known results. In particular, after the appearance of the articles \cite{AAP2016} and  \cite{KS2017}, several investigations and
new problems on Bohr's inequality in the plane case appeared in the literature (cf. \cite{AAL20,BhowDas-18,LP2018,KayPon3, LSX2018,PVW2019,PW2019}).

One of our aims in this article is to address the harmonic analog of this question (see Problem \ref{LSH-prob1})
raised by  Paulsen et al. \cite{PPS2002} but with a refined formulation as in \cite{PW2019} (see Theorem \Ref{Theo-ABC}).

\subsection{Classical Inequality of H. Bohr}
Let $\mathcal{B}$ be the Schur class of all analytic functions $f$ on the open unit disk
$\mathbb{D}:=\{z \in \mathbb{C}:\,|z|<1\}$ such that $\|f\|_\infty :=\sup_{z\in \ID}|f(z)|\leq 1$.
Then the classical inequality
examined by Bohr in 1914 \cite{B1914} states that $1/3$ is the largest value of $r\in [0,1)$ for which
the following inequality holds:
\begin{equation}
B(f,r):=\sum_{k=0}^{\infty}|a_k|r^k\leq 1
\label{liu1}
\end{equation}
for every analytic function $f\in \mathcal{B}$ with the Taylor series expansion $f(z)=\sum_{k=0}^{\infty}a_kz^k$.
Bohr actually obtained that  \eqref{liu1} is true when  $r\leq 1/{6}$. Later   Riesz, Schur and Wiener, independently
established the Bohr inequality \eqref{liu1} for $r\leq 1/{3}$ and that $1/3$
is the best possible constant. It is quite natural that the constant $1/3$ is called the Bohr radius
for the space $\mathcal{B}$. Moreover, for
$$\varphi_a(z)=\frac{a-z}{1-a z},\quad a\in [0,1),
$$
it follows easily that $B(\varphi_a,r)>1$ if and only if $r>1/(1+2a)$, which for $a\ra 1$ shows that $1/3$ is optimal.
Bohr's and Wiener's proofs can be found in \cite{B1914}. Other proofs of Bohr's inequality may be
found from \cite{S1927,T1962}.  Then it is worth pointing out that there is no extremal function in $\mathcal{B}$ such that the Bohr radius
is precisely $1/3$ (cf. \cite[Corollary 8.26]{GarMasRoss-2018}).

\subsection{The Bohr radius for functions having multiple zeros at the origin}
\bprob\label{LSH-prob1}
In \cite{PPS2002}, the authors considered among others for $k\in\mathbb{N}$ the classes
$\mathcal{B}_{k}:=z^k\mathcal{B}$, that is,
$$
\mathcal{B}_{k}=\left\{f \in \mathcal{B}: \,f(0)=\cdots=f^{(k-1)}(0)=0\right\} :=\left\{z^kf:\, f \in \mathcal{B}\right\},
$$
and asked for which $r_k\in(0, 1)$, and
\begin{equation}
f(z)=\sum_{n=k}^{\infty} a_{n} z^{n} 
\in \mathcal{B}_{k}
\label{liu2}
\end{equation}
we have  the inequality
\begin{equation}
\sum_{n=k}^{\infty}|a_n|r^n\leq 1 ~\mbox{ for $r\in [0,r_k]$}
\label{liu3}
\end{equation}
and for each $r\in (r_k,1)$  there exists a function $f_k(z)=\sum_{n=k}^{\infty} a_{n}^{(k)} z^{n}$ in $\mathcal{B}_{k}$
such that $\sum_{n=k}^{\infty}|a_{n}^{(k)}|r^n > 1$. Here the constant $r_k$ is referred to as the Bohr radius of order $k$.
\eprob

Clearly, $\mathcal{B}_{0}=\mathcal{B}$, and $\mathcal{B}_{1}=\left\{f \in \mathcal{B}: \,f(0)=0\right\}$. For
$f\in\mathcal{B}_{1}$ (i.e. for $k=1)$, Tomi\'{c} \cite{T1962} proved that (\ref{liu3}) holds for $0 \leq r \leq 1 / 2$
(also obtained by Landau independently, see \cite{LD1986}). Later Ricci \cite{RC1955} established that
this holds for $0 \leq r \leq 3 / 5$, and the largest value of $r$ for which (\ref{liu3}) holds would lie
in the interval $(3/5, 1/\sqrt{2}]$. Later in 1962,
Bombieri \cite{BS1962} found that the inequality (\ref{liu3}) holds for $r \in[0,1 / \sqrt{2}]$, where the upper
bound cannot be improved. An alternate proof of this result may be found from a recent paper of Kayumov and Ponnusamy
\cite{KayPon_AAA18} in which they solved an open problem of Djakov and Ramanujan on powered Bohr inequality. However,
Problem \ref{LSH-prob1} for $k\geq 2$ remains open. On the other hand, in connection with Problem \ref{LSH-prob1},
Ponnusamy and Wirths \cite{PW2019} proved the following sharp inequalities for $k\geq 2$ while the case $k=1$
has been proved in \cite{PVW2019}:

\begin{Thm}\label{Theo-ABC}
For $k\geq 1$, let $f\in \mathcal{B}_{k}$ have an expansion \eqref{liu2} and
$$M_k^{1}(f,r)=\sum_{n=k}^{\infty}\left|a_{n}\right| r^{n}+\left(\frac{1}{1+\left|a_{k}\right|}+\frac{r}{1-r}\right) \sum_{n=k+1}^{\infty}\left|a_{n}\right|^{2} r^{2 n -k}
$$
and
$$M_k(f,r)=\sum_{n=k}^{\infty}\left|a_{n}\right| r^{n}+\left(\frac{1}{1+\left|a_{k}\right|}+\frac{r}{1-r}\right) \sum_{n=k}^{\infty}\left|a_{n}\right|^{2} r^{2 n-k}.
$$
Then we have the following inequalities:
\begin{enumerate}
\item $M_k^{1}(f,r) \leq 1 $
is valid for $r \in\left[0, R_{k}\right],$ where $R_{k}$ is the unique root in $(0,1)$ of the equation
\begin{equation*}
4(1-r)-r^{k-1}\left(1-2 r+5 r^{2}\right)=0.
\end{equation*}
The upper bound $R_{k}$ cannot be improved.
\item $M_k(f,r) \leq 1 $ is valid for $r \in\left[0, S_{k}\right]$, where $S_{k}$ is the unique root in $(0,1)$ of the equation
\begin{equation*}
2(1-r)-r^{k}(3-r)=0.
\end{equation*}
The upper bound $S_{k}$ cannot be improved. Also, as $M_k^{1}(f,r)\leq M_k(f,r)$, it follows that $S_k\leq R_k$.

\item With $\left|a_{k}\right|=a \in(0,1]$ being fixed, $M_k(f,r) \leq 1 $
is valid for $r \in\left[0, \rho_{k}(a)\right]$, where $\rho_{k}(a)$ is the unique root in $(0,1)$ of the equation
\begin{equation*}
(1+a)(1-r)-r^{k}\left[2 a^{2}+a+r\left(1-2 a^{2}\right)\right]=0.
\end{equation*}
The upper bound $\rho_{k}(a)$ cannot be improved.
\end{enumerate}
\end{Thm}

\br
We note that $\rho_{k}(1)=S_k$. As $M_k^{1}(f,r)\leq M_k(f,r)$, it follows that $S_{k}\leq R_{k}$.
\er

\subsection{The Bohr radius for $p$-symmetric functions}

Recently, Kayumov  et al. \cite{KS2017} have obtained the following general result. As a corollary to this,
an open problem raised by Ali et al. \cite{AliBarSoly} about the determination of Bohr radius for odd functions
from $\mathcal{B}$ has been settled affirmatively.

\begin{Thm}\label{Theo-B}
(\cite{KS2017})
Let $m, \,p\in \mathbb{N}$,  $ m\leq p$, and $f\in \mathcal{B}$ with $f(z)=\sum_{k=0}^{\infty}a_{pk+m}z^{pk+m}$.
Then
\begin{equation*}
B(f,r)= \sum_{k=0}^{\infty}|a_{pk+m}|r^{pk+m} \leq1~ \mbox{ for $r\leq r_{p,m}$},
\end{equation*}
where $r_{p,m}$ is the maximal positive root of the equation
$ -6r^{p-m}+r^{2(p-m)}+8r^{2p}+1=0.
$
The extremal function has the form $z^{m}(z^{p}-a)/(1-az^{p})$, where
\begin{equation*}
a=\left (1-\frac{\sqrt{1-r_{p,m}^{2p}}}{\sqrt{2}}\right )\frac{1}{r_{p,m}^{p}}.
\end{equation*}
\end{Thm}

\br
We note that the case $m=0$ is trivial as it follows from the classical theorem of H. Bohr with a change of variable $\zeta =z^p$.
This gives the condition $r\leq r_{p,0}=1/\sqrt[p]{3}$. The case $p=2$ and $m=1$ corresponds to the question raised
by Ali et al. \cite{AliBarSoly}.
\er

\subsection{The Bohr radius for harmonic functions}

In \cite{KSS2017}, the authors initiated the discussion on Bohr radius for the class of
complex-valued function $f = u + iv$ harmonic in $\ID$, where $u$ and $v$ are real-valued harmonic functions of $\ID$.
It follows that $f$ admits the canonical representation $f=h+\overline{g}$,
where $h$ and $g$ are analytic in $\ID$  such that $f(0)=0=g(0)$.
The Jacobian $J_f (z)$ of $f$ is given by $J_f (z)=|h'(z)|^2-|g'(z)|^2$, and we say that a
locally univalent harmonic function $f$ in $\mathbb{D}$ is said to be sense-preserving if  $J_f (z)>0$ in $\ID$;
or equivalently, its dilatation $\omega =g'/h'$ is an analytic function in $\ID$
which maps  $\mathbb{D}$ into itself (cf. \cite{D1983} or \cite{L1936}).

If a locally univalent and sense-preserving harmonic mapping $f=h+\overline{g}$ satisfies the condition
$\left|\omega (z)\right|\leq d<1 $ in $\ID$,
then $f$ is called   $K$-quasiregular harmonic mapping on $\mathbb{D}$, where $K=\frac{1+d}{1-d}\geq 1$
(cf. \cite{Kalaj2008,Martio1968}). Obviously $d\rightarrow 1$ corresponds to the case $K\rightarrow \infty$.

For a harmonic function $f=h+\overline{g}$ in $\ID$, where $h$ and $g$ admit power series
expansions of the form $h(z) =\sum_{n=0}^{\infty}a_n z^n$ and $g(z)=\sum_{n=0}^{\infty}b_n z^n$,
we denote the classical Bohr sum by
$$ B_H(f,r) := B(h,r) + B(g,r) =\sum_{n=0}^{\infty}(|a_{n}|+|b_{n}|)r^{n}.
$$
A harmonic function $f=h+\overline{g}$ in $\ID$ is said to be $p$-symmetric if $h$ and $g$ have the form
$h(z)=\sum_{n=0}^{\infty} a_{n}z^{pn+m}$ and $g(z)=\sum_{n=0}^{\infty} b_{n}z^{pn+m}$
for some $m\in \IN_0=\IN\cup \{0\}$.
Harmonic extension of the classical Bohr theorem was established first in \cite{KSS2017}. For example, they proved the
following result (Theorem \Ref{Theo-C}). Furthermore, the Bohr radii for harmonic  and starlike log harmonic mappings
in $\ID$ were investigated, for example, in \cite{EPR-2017,KS2017,KSS2017,LP2019}, and in some cases in improved form.

\begin{Thm}\label{Theo-C}
(\cite{KSS2017})
 Let $p\in \IN$ and $p\geq2$. Suppose that $f(z)=h(z)+\overline{g(z)}=\sum_{n=0}^{\infty} a_{n}z^{pn+1}+\sum_{n=0}^{\infty }\overline{b_{n}z^{pn+1}}$
 is a harmonic $p$-symmetric function in $\ID$, where $h$ and $g$ are bounded functions in $\ID$. Then
\begin{equation*}
 B_H(f,r) =\sum_{n=0}^{\infty}(|a_{n}|+|b_{n}|)r^{pn+1}\leq\max\{\|h\|_{\infty},\|g\|_{\infty}\}~\mbox{ for $r\leq 1/2$}.
\end{equation*}
The number $1/{2}$ is sharp.
\end{Thm}

It is natural to raise the following.

\bprob\label{HLP-prob1}
Whether Theorem~C 
 holds under a weaker hypotheses, namely, by replacing the condition ``boundedness of
 $h$ and $g$" by ``$|g'(z)| \leq |h'(z)|$ and $h$ is bounded."
\eprob

In Theorem \ref{HLP-th2-new}, 
we present an affirmative answer to this question in a more general setting.

The paper is organized as follows.  In Section \ref{HLP-sec2}, we present the main results of this paper.
In Theorem \ref{HLP-th2-new}, we present an affirmative answer to Problem \ref{HLP-prob1} in a general form, and Corollary
\ref{HLP-th2} answers Problem \ref{HLP-prob1}.  As consequence,  generalization Theorem~C 
(with of higher order zero at the origin) is established (see Theorem \ref{HLP-th1}).
In Section \ref{HLP-sec3}, we  state and prove several lemmas. In addition, we present the proof of Bohr's
inequalities for the class of harmonic mappings, which improve
the first two items in Theorems~A and C.

In Section \ref{HLP-sec4}, we state and prove three theorems which extend three recent results of
Ponnusamy et al. \cite{PW2019} from the case of analytic functions to the case of sense-preserving harmonic mappings.

\section{Main Results}\label{HLP-sec2}

We now state a generalization of Theorem~C 
in a general setting and the next result (Theorem \ref{HLP-th1}) is a direct generalization of Theorem~C.

\bthm\label{HLP-th2-new}
Let $m, \, p\in \mathbb{N}$,  $p\geq2$. Suppose that $f(z)=h(z)+\overline{g(z)}=\sum_{k=0}^{\infty} a_{k}z^{pk+m}+\sum_{k=0}^{\infty }\overline{b_{k}z^{pk+m}}$ is harmonic and $p$-symmetric in $\ID$ such that $|h^{(m)}(0)|=|g^{(m)}(0)|$ and $|g'(z)|\leq d|h'(z)|$ in $\ID\backslash \{0\}$ for some $d\in [0,1]$, where $h$ is bounded. Then the following hold:

\begin{enumerate}
\item If $\frac{p}{m}>\log_2 (2+d)$, then
\begin{equation*}
  B_H(f,r)=\sum_{k=0}^{\infty}(|a_{k}|+|b_{k}|)r^{pk+m}\leq\|h\|_{\infty}~ \mbox{ for  } ~r\leq\sqrt[m]{\frac{1}{2}}.
\end{equation*}
When $d=1$, the extremal mapping has the form $f(z)=h(z)+\overline{\lambda h(z)}$ with $h(z)=z^m$ and $|\lambda|=1$.

\item If $1\leq\frac{p}{m}\leq\log_2 (2+d)$, then
\begin{equation*}
B_H(f,r)\leq\|h\|_{\infty}~\mbox{ for  }~r\leq r_{p,m,d},
\end{equation*}
where $r_{p,m,d}$ is the maximal positive root of the equation
\begin{equation}
  r^{2(p-m)}-(8+4d)r^{p-m} +4(1+d)(3+d)r^{2p}+4=0.
  \label{liu31}
\end{equation}
When $d=1$, the extremal function is given by $f(z)=h(z)+\overline{\lambda h(z)}$,  $|\lambda|=1$, where
$$h(z)=z^{m}\left (\frac{z^{p}-a}{1-az^{p}}\right ),~\mbox{ with }~
a=\left (1-
\frac{\sqrt{1-r_{p,m,1}^{2p}}}{\sqrt{2}}\right )\frac{1}{r_{p,m,1}^{p}}.
$$
\end{enumerate}
\ethm

\bcor\label{HLP-cor2} Suppose that $m, \, p\in \mathbb{N}$, and $f(z)=\sum_{k=0}^{\infty}a_{pk+m}z^{pk+m} \in\mathcal{B}_{pk+m}$.
\begin{enumerate}
\item If $1\leq \frac{p}{m}\leq\log_{2}3 \approx 1.58496$, then
\begin{equation*}
  B(f,r)\leq\frac{1}{2}~\mbox{ for }~ r\leq r_{p,m},
\end{equation*}
where $r_{p,m}$ is the maximal positive root of the equation
\begin{equation}\label{liu01}
  -12r^{p-m}+r^{2(p-m)}+32r^{2p}+4=0.
\end{equation}
 The extremal function is given by
\be\label{HLP-eq2}
f(z)=z^{m}\left (\frac{z^{p}-a}{1-az^{p}}\right ),~\mbox{ with }~a=\left (1-\frac{\sqrt{1-r_{p,m}^{2p}}}{\sqrt{2}}\right )\frac{1}{r_{p,m}^{p}}.
\ee

\item  If $\frac{p}{m}>\log_{2}3$, then
\begin{equation*}
  B(f,r)\leq\frac{1}{2}~~\mbox{ for  }~~r\leq \sqrt[m]{\frac{1}{2}}.
\end{equation*}
The extremal function has the form $z^{m}$.
\end{enumerate}
\ecor
\bpf
Apply the method of the proof of Theorem \ref{HLP-th2-new} (by setting $d=1$, $g(z)\equiv 0$).
\epf

We now state a direct generalization of Theorem~C. 

\bthm\label{HLP-th1}
Let $m, \, p\in \mathbb{N}$,  $p\geq2$. Suppose that
$f(z)=h(z)+\overline{g(z)}=\sum_{k=0}^{\infty} a_{pk+m}z^{pk+m}$ $+\sum_{k=0}^{\infty }\overline{b_{pk+m}z^{pk+m}}$
is harmonic in $\ID$, where $h$ and $g$ are bounded. The following hold:

\begin{enumerate}
\item If $\frac{p}{m}>\log_{2}3\approx 1.58496$, then
$$B_H(f,r) \leq \max\{\|h\|_{\infty},\|g\|_{\infty}\}~ \mbox{ for $\ds r\leq \sqrt[m]{\frac{1}{2}}$}.
$$
The extremal function is given by $f(z)=z^{m} +\overline{\lambda z^{m}}$, $|\lambda|=1$.

\item If $1\leq \frac{p}{m}\leq\log_{2}3$, then
$$ B_H(f,r)  \leq \max\{\|h\|_{\infty},\|g\|_{\infty}\} ~\mbox{ for  $r\leq r_{p,m}$},
$$
where $r_{p,m}$ is the maximal positive root in $(0,1)$ of the equation  \eqref{liu01}.

The extremal function is given by $f(z)=h(z) +\overline{\lambda h(z)}$, $|\lambda|=1$, where
$$h(z)=z^{m}\left (\frac{z^{p}-a}{1-az^{p}}\right ),~\mbox{ with }~a=\left (1-\frac{\sqrt{1-r_{p,m}^{2p}}}{\sqrt{2}}\right )\frac{1}{r_{p,m}^{p}}.
$$
\end{enumerate}
\ethm

\br
If we set $m=1$ in Theorem \ref{HLP-th1}(1), then we get Theorem~C. 
\er

Note that the following corollary generalizes Theorem \ref{HLP-th1} under the conditions ``$|h'(0)|=|g'(0)|$ and $|g'(z)|\leq|h'(z)|$ in
$\ID\backslash \{0\}$" instead of ``$h$ and $g$ being bounded in $\ID$."

\bcor\label{HLP-th2}
Let $m, \, p\in \mathbb{N}$,  $p\geq2$. Suppose that $f(z)=h(z)+\overline{g(z)}=\sum_{k=0}^{\infty} a_{k}z^{pk+m}+\sum_{k=0}^{\infty }\overline{b_{k}z^{pk+m}}$ is harmonic and $p$-symmetric in $\ID$ such that $|h'(0)|=|g'(0)|$ and $|g'(z)|\leq|h'(z)|$ in $\ID\backslash \{0\}$,
where $h$ is bounded. Then the conclusions {\rm (1)} and  {\rm (2)} of Theorem \ref{HLP-th1} continue to hold.
\ecor
\bpf
Set $d=1$ in Theorem \ref{HLP-th2-new} and let $ r_{p,m}:=r_{p,m,1}$.
\epf

Because of its independent interest, let us next state the following result as a corollary to  Theorems~A 
Indeed, applying the analogous methods as in the proofs of the three cases of Theorems~A, 
we have the following. So we omit the details.

\bcor\label{HLP-cor3}
For $k\geq 1$, $m, \,p\in \mathbb{N}$ and $m \leq p$,  we let $f(z)=\sum_{n=k}^{\infty} a_{p n+m} z^{p n+m} \in\mathcal{B}_{pk+m}$ and
$$M_{pk+m}^{1}(f,r)=
\sum_{n=k}^{\infty}\left|a_{p n+m}\right| r^{p n+m}+\left(\frac{1}{1+\left|a_{p k+m}\right|}+\frac{r^{p}}{1-r^{p}}\right) \sum_{n=k+1}^{\infty}\left|a_{p n+m}\right|^{2} r^{p(2 n-k)+m}
$$
and
$$M_{pk+m}(f,r)=\sum_{n=k}^{\infty}\left|a_{p n+m}\right| r^{p n+m}+\left(\frac{1}{1+\left|a_{p k+m}\right|}+\frac{r^{p}}{1-r^{p}}\right) \sum_{n=k}^{\infty}\left|a_{p n+m}\right|^{2} r^{p(2 n-k)+m}.
$$
 Then we have the following inequalities:
\begin{enumerate}
\item $M_{pk+m}^{1}(f,r)\leq 1$ is valid for $r \in\left[0, v_{k}\right]$, where $v_{k}$ is the unique root in $(0, 1)$ of the equation
\begin{equation*}
4\left(1-r^{p}\right)-r^{p(k-1)+m}\left(5 r^{2 p}-2 r^{p}+1\right)=0.
\end{equation*}
The upper bound $v_{k}$ cannot be improved.
\item $M_{pk+m}(f,r)\leq 1$
is valid for $r \in\left[0, \omega_{k}\right]$, where $\omega_{k}$ is the unique root in $(0, 1)$ of the equation
\begin{equation*}
2\left(1-r^{p}\right)-r^{p k+m}\left(3-r^{p}\right)=0.
\end{equation*}
The upper bound $\omega_{k}$ cannot be improved.

\item With $\left|a_{p k+m}\right|=a \in(0,1]$ being fixed, $M_{pk+m}(f,r)\leq 1$
is valid for $r \in\left[0, \, \eta_{k}\right]$, where $\eta_{k}$ is the unique root in $(0,1)$ of the equation
\begin{equation*}
(1+a)\left(1-r^{p}\right)-r^{p k+m}\left[2 a^{2}+a+r^{p}\left(1-2 a^{2}\right)\right]=0.
\end{equation*}
The upper bound $\eta_{k}$ cannot be improved.
\end{enumerate}
\ecor
\bpf
The desired conclusion follows if we write  $f(z)$ as $f(z)=z^mt(z^p)$,
where $t(z)=\sum_{n=k}^{\infty} a_{p n+m} z^n \in\mathcal{B}_k$,
and apply the proof of Theorem~A. 
\epf

Next we generalize Theorem~A 
or Corollary \ref{HLP-cor3} by establishing Bohr-type inequalities for harmonic mappings
with multiple zero at the origin.

\bthm\label{HLP-th3}
Let $k\geq 1$, $m, \, p\in \mathbb{N}$ and $m\leq p$. Suppose that $f=h+\overline{g}$ is harmonic in $\ID$, where $h$ and $g$
are given by
\begin{eqnarray}
h(z)=\sum_{n=k}^{\infty} a_{pn+m}z^{pn+m} ~\mbox{ and }~g(z)=\sum_{n=k}^{\infty }b_{pn+m}z^{pn+m}.
\label{liu50}
\end{eqnarray}
In addition, let $|g'(z)|\leq d|h'(z)|$ in $\ID\backslash \{0\}$ for some $d\in [0,1]$ and $h\in\mathcal{B}_{pk+m}$.
Define
\be\label{HLP-eq6}
M_{pk+m}(h,r)=\sum_{n=k}^{\infty}\left|a_{p n+m}\right| r^{p n+m}+\left(\frac{1}{1+\left|a_{p k+m}\right|}+\frac{r^{p}}{1-r^{p}}\right) \sum_{n=k}^{\infty}\left|a_{p n+m}\right|^{2} r^{p(2 n-k)+m},
\ee
and
$$N_{pk+m}(g,r)=\sum_{n=k}^{\infty}\left|b_{p n+m}\right| r^{p n+m}+\left(\frac{1}{1+\left|a_{p k+m}\right|}+\frac{r^{p}}{1-r^{p}}\right) \sum_{n=k}^{\infty}\left|b_{p n+m}\right|^{2} r^{p(2 n-k)+m}.
$$
Then the inequality
\begin{eqnarray}
M_{pk+m}(h,r)+ N_{pk+m}(g,r)\leq 1
\label{liu501}
\end{eqnarray}
is valid for $r \in\left[0, r_k\right]$, where $r_k=\min\{r_{k}',  1/\sqrt[p]{3}\}$, and $r_{k}'$ is the unique root in $(0, 1)$ of the equation
$t_{k}(r)=0$, where
\be\label{HLP-eq5}
t_{k}(r)=\frac{2}{d+1}(1-r^{p})-r^{p k+m}\left(3-r^{p}\right).
\ee
\ethm

\bthm\label{HLP-th4}
Let $k\geq 2$, $m, \, p\in \mathbb{N}$ and $m\leq p$. Suppose that $f=h+\overline{g}$ is harmonic in $\ID$, where $h$ and $g$
are given by  \eqref{liu50},
In addition, let $h,\, g\in\mathcal{B}_{pk+m}$ and define
$$M_{pk+m}^{1}(h,r)=\sum_{n=k}^{\infty}\left|a_{p n+m}\right| r^{p n+m}+\left(\frac{1}{1+\left|a_{p k+m}\right|}+\frac{r^{p}}{1-r^{p}}\right) \sum_{n=k+1}^{\infty}\left|a_{p n+m}\right|^{2} r^{p(2 n-k)+m}.
$$
Then the inequality
\be\label{liu51}
M_{pk+m}^{1}(h,r)+ M_{pk+m}^{1}(g,r)\leq 1
\ee
is valid for $r \in\left[0, \tau_{k}\right]$, where $\tau_{k}$ is the unique root in $(0,1)$ of the equation
\begin{equation*}
2\left(1-r^{p}\right)-r^{p(k-1)+m}\left(5 r^{2 p}-2 r^{p}+1\right)=0.
\end{equation*}
The upper bound $\tau_{k}$ cannot be improved.
\ethm

\bthm\label{HLP-th5}
Let $k\geq 2$, $m, \, p\in \mathbb{N}$ and $m\leq p$. Suppose that $f=h+\overline{g}$ is harmonic in $\ID$, where $h$ and $g$
are given by  \eqref{liu50}, and $h,\, g\in\mathcal{B}_{pk+m}$. Then the inequality
\begin{eqnarray}
M_{pk+m}(h,r)+ M_{pk+m}(g,r)\leq 1
\label{liu53}
\end{eqnarray}
is valid for $r \in\left[0, \theta_{k}\right]$, where $M_{pk+m}(h,r)$ is given by \eqref{HLP-eq6} and
$\theta_{k}$ is the unique root in $(0,1)$ of the equation
\begin{equation*}
\left(1-r^{p}\right)-r^{p k+m}\left(3-r^{p}\right)=0.
\end{equation*}
The upper bound $\theta_{k}$ cannot be improved as $f(z)=h(z) +\overline{\lambda h(z)}$ shows,
where $h(z)=z^{p k+m}$ and $|\lambda|=1$.
\ethm

Our next result is similar to Theorem \ref{HLP-th4}, but for fixed initial coefficients $a_{p k+m}$ and $b_{p k+m}$ having
same modulus value.

\bthm \label{HLP-th6}
Let $k\geq 2$, $m, \, p\in \mathbb{N}$ and $m\leq p$. Suppose that $f=h+\overline{g}$ is harmonic in $\ID$, where $h$ and $g$
are given by  \eqref{liu50}, and $h,\, g\in\mathcal{B}_{pk+m}$. Let $\left|a_{p k+m}\right|=\left|b_{p k+m}\right|=a \in(0,1]$
be fixed. Then the inequality
\begin{eqnarray}
M_{pk+m}(h,r)+ M_{pk+m}(g,r)\leq 1
\label{liu55}
\end{eqnarray}
is valid for $r \in\left[0, \varsigma_{k}\right]$, where $\varsigma_{k}=\varsigma_{k}(a)$ is the unique root in $(0,1)$ of the equation
\begin{equation*}
(1+a)\left(1-r^{p}\right)-2r^{p k+m}\left[2 a^{2}+a+r^{p}\left(1-2 a^{2}\right)\right]=0.
\end{equation*}
The upper bound $\varsigma_{k}$ cannot be improved.
\ethm

\br
Clearly, $\varsigma_{k}(1)=\theta_{k}$. Also, it is possible to fix both $\left|a_{p k+m}\right|$ and $\left|b_{p k+m}\right|$
separately and obtain an analogous general result than Theorem \ref{HLP-th6}.
\er

\section{Key lemmas and their Proofs}\label{HLP-sec3}
In order to establish our main results, we need the following lemmas. 

\blem\label{HLP-lem1} 
Suppose that $m, \, p\in \mathbb{N}$, $m\leq p$, $d\in (0, 1]$, and $r=r_{p,m,d}$ is as in Theorem \ref{HLP-th2-new}, i.e. the root of \eqref{liu31} in $(0,1)$. Then
$$ r^{p+m}\leq\frac{1}{3+d}.
$$
\elem
\bpf
 Let $y=r^{p+m}_{p,m,d}$. Then \eqref{liu31} becomes a quadratic equation in $y$ of the form
\begin{equation*}
  \left (4(1+d)(3+d)+\frac{1}{r^{2m}_{p,m,d}}\right )y^{2}-(8+4d)y+4r^{2m}_{p,m,d}=0,
\end{equation*}
which has two solutions
\begin{eqnarray*}
y&=&\frac{4+2d\pm 2\sqrt{(1+d)(3+d)}\sqrt{1-4r^{2m}_{p,m,d}}}{4(1+d)(3+d)+\frac{1}{r^{2m}_{p,m,d}}}\\
&\leq &\frac{4+2d+ 2\sqrt{(1+d)(3+d)}\sqrt{1-4r^{2m}_{p,m,d}}}{4(1+d)(3+d)+\frac{1}{r^{2m}_{p,m,d}}}\\
&\leq& \frac{1}{2}\left (\sup \frac{2+d+ \sqrt{(1+d)(3+d)}\sqrt{1-4r^{2m}_{p,m,d}}}{(1+d)(3+d)+\frac{1}{4r^{2m}_{p,m,d}}}\right )
\end{eqnarray*}
and therefore, it is a simple exercise to see that
\begin{eqnarray*}
r^{p+m}_{p,m,d}&\leq& \frac{1}{2}\left (\sup_{t\in(0,1]} \frac{2+d+ \sqrt{(1+d)(3+d)}\sqrt{1-t}}{(1+d)(3+d)+\frac{1}{t}}\right )\\
&=& \left . \frac{1}{2}\left (\frac{2+d+ \sqrt{(1+d)(3+d)}\sqrt{1-t}}{(1+d)(3+d)+\frac{1}{t}}\right )\right |_{t=\frac{2}{3+d}}\\
& = & 
\frac{1}{3+d},
\end{eqnarray*}
which completes the proof of the lemma. 
\epf

\blem\label{HLP-lem2} 
Suppose that $m, \, p\in \mathbb{N}$, $m\leq p$, $d\in (0, 1]$, and $r=r_{p,m,d}$ is as in Theorem \ref{HLP-th2-new},
i.e. the root of \eqref{liu31} in $(0,1)$. Then
\begin{equation*}
 \frac{1}{r^{p-m}}(2+d-\sqrt{(1+d)(3+d)}\sqrt{1-r^{2p}})=\frac{1}{2}.
\end{equation*}
\elem
\bpf
Suppose that $m<p$ and let $y=r^{p-m}$. Then \eqref{liu31} reduces to a quadratic equation in $y$
\begin{equation*}
  y^2-(8+4d)y+4(1+d)(3+d)r^{2p}+4=0,
\end{equation*}
which has two solutions
$$y_{1}=4+2d+2\sqrt{(1+d)(3+d)}\sqrt{1-r^{2p}}>1~\mbox{ and }~
y_{2}=4+2d-2\sqrt{(1+d)(3+d)}\sqrt{1-r^{2p}}.
$$
The solution $y=y_{1}$ is impossible because all positive roots of the initial equation must be less than $1$. Therefore,
$$y=y_{2}=2\left (2+d-\sqrt{(1+d)(3+d)}\sqrt{1-r^{2p}}\right ).
$$

Now, consider the case $m=p$. In this case
$$r_{m,m,d}=\left (\frac{3+4d}{4(1+d)(3+d)}\right )^{\frac{1}{2m}}
$$
so that

\vspace{6pt}

$\ds \frac{1}{r^{p-m}}\left (2+d-\sqrt{(1+d)(3+d)}\sqrt{1-r^{2m}}\right)$
\begin{eqnarray*}
&=&2+d-\sqrt{(1+d)(3+d)}\sqrt{1-\frac{3+4d}{4(1+d)(3+d)}}\\
&=&
2+d-\left (\frac{2d+3}{2}\right )=\frac{1}{2},
\end{eqnarray*}
and the proof is complete.
\epf

\blem\label{HLP-lem4} (\cite{KP2017})
Let $0<R\leq1$. If $g(z)=\sum_{k=0}^{\infty}b_{k}z^{k}$ is analytic and satisfies the inequality $|g(z)|\leq1$ in $\ID$. Then the following sharp inequality holds:
\begin{equation}
  \sum_{k=1}^{\infty}|b_{k}|^{2}R^{pk}\leq R^{p}\frac{(1-|b_{0}|^{2})^{2}}{1-|b_{0}|^{2}R^{p}}.
  \label{liu02}
\end{equation}
\elem

\blem\label{HLP-lem6} (\cite{PVW2019})
If $f \in\mathcal{B}$ has the expansion $f(z)=\sum_{n=0}^{\infty} a_{n} z^{n}$, then
$$\sum_{n=0}^{\infty}\left|a_{n}\right| r^{n}+\left(\frac{1}{1+\left|a_{0}\right|}+\frac{r}{1-r}\right) \sum_{n=1}^{\infty}\left|a_{n}\right|^{2} r^{2 n} \leq |a_0|+\frac{r}{1-r}(1-|a_0|^2).$$
\elem

\section{Bohr's inequality for the class of harmonic mappings}\label{HLP-sec4}

\subsection{Proof of Theorem \ref{HLP-th2-new}}
Given that $|g'(z)|\leq d|h'(z)|$ for some $d\in (0,1]$, where
$$h(z)= \sum_{k=0}^{\infty} a_{k}z^{pk+m} ~\mbox{ and }~ g(z)=\sum_{k=0}^{\infty }b_{k}z^{pk+m}.
$$
We integrate inequality $|g'(z)|^{2}\leq d^2|h'(z)|^{2}$ over the circle $|z|=r$ and get
\begin{equation*}
  \sum_{k=0}^{\infty}(pk+m)^{2}|b_{k}|^{2}r^{2(pk+m-1)}\leq d^2\sum_{k=0}^{\infty}(pk+m)^{2}|a_{k}|^{2}r^{2(pk+m-1)}.
\end{equation*}
We integrate the last inequality with respect to $r^{2}$ and obtain
\begin{equation*}
  \sum_{k=0}^{\infty}(pk+m)|b_{k}|^{2}r^{2(pk+m)}\leq d^2\sum_{k=0}^{\infty}(pk+m)|a_{k}|^{2}r^{2(pk+m)}.
\end{equation*}
One more integration (after dividing by $r^{2}$) gives
\begin{equation*}
  \sum_{k=0}^{\infty}|b_{k}|^{2}r^{2(pk+m)}\leq d^2\sum_{k=0}^{\infty}|a_{k}|^{2}r^{2(pk+m)},
\end{equation*}
which (since $|a_0|=|b_0|$ by hypothesis) yields
\begin{equation}
  \sum_{k=1}^{\infty}|b_{k}|^{2}r^{pk}\leq d^2\sum_{k=1}^{\infty}|a_{k}|^{2}r^{pk} ~\mbox{ for $r<1$.}
  \label{liu08}
\end{equation}

For simplicity, we suppose that $\|h\|_{\infty}=1$.

Following the idea from \cite{KP2017} (see also \cite[Proof of Theorem 1]{KS2017}), one can obtain firstly that
\begin{eqnarray} \label{liu03-e}
B(h,r)=r^{m}\sum_{k=1}^{\infty}|a_{k}|r^{pk} &\leq& \frac{r^{p+m}(1-a^{2})}{\sqrt{1-a^{2}r^{p}\rho^{p}}}\frac{1}{\sqrt{1-\rho^{-p}r^{p}}},
\end{eqnarray}
where $a=|a_0|$, and for any $\rho>1$ such that $\rho r\leq1$. Indeed, we may let $h(z)=z^{m}t(z^{p})$,
where $t(z)=\sum_{k=0}^{\infty}a_{k}z^{k}\in\mathcal{B}$.
Also, let  $r=r_{p,m,d}$ and  $|a_{0}|=a$. Then, 
as in  \cite[Proof of Theorem 1]{KS2017}, it follows easily that
\begin{eqnarray}\label{HLP-eq3}
\sum_{k=1}^{\infty}|a_{k}|r^{pk} &\leq& \sqrt{\sum_{k=1}^{\infty}|a_{k}|^{2}\rho^{pk}r^{pk}}\sqrt{\sum_{k=1}^{\infty}\rho^{-pk}r^{pk}} \nonumber\\
&\leq& \sqrt{r^{p}\rho^{p}\frac{(1-a^{2})^{2}}{1-a^{2}r^{p}\rho^{p}}}\sqrt{\frac{\rho^{-p}r^{p}}{1-\rho^{-p}r^{p}}}
\, =\frac{r^{p}(1-a^{2})}{\sqrt{1-a^{2}r^{p}\rho^{p}}}\frac{1}{\sqrt{1-\rho^{-p}r^{p}}},
\end{eqnarray}
for any $\rho>1$ such that $\rho r\leq1$.
In the first and the second steps above we have used the classical Cauchy-Schwarz inequality, and (\ref{liu02}) with $R=\rho r$ in
Lemma \ref{HLP-lem4}, respectively. Hence, \eqref{liu03-e} follows.

Secondly, using the classical Cauchy-Schwarz inequality and \eqref{liu08},
we see that
\begin{eqnarray*}
\sum_{k=1}^{\infty}|b_{k}|r^{pk} &\leq& \sqrt{\sum_{k=1}^{\infty}|b_{k}|^{2}\rho^{pk}r^{pk}}\sqrt{\sum_{k=1}^{\infty}\rho^{-pk}r^{pk}}\\
&\leq& d \sqrt{\sum_{k=1}^{\infty}|a_{k}|^{2}\rho^{pk}r^{pk}}\sqrt{\frac{\rho^{-p}r^{p}}{1-\rho^{-p}r^{p}}}  ~\mbox{ (by \eqref{liu08})}
\end{eqnarray*}
and thus, by \eqref{HLP-eq3} we have
\begin{eqnarray} \label{liu03-f}
B(g,r)=r^{m}\sum_{k=1}^{\infty}|b_{k}|r^{pk}
 &\leq& \frac{dr^{p+m}(1-a^{2})}{\sqrt{1-a^{2}r^{p}\rho^{p}}}\frac{1}{\sqrt{1-\rho^{-p}r^{p}}}
\end{eqnarray}
for any $\rho>1$ such that $\rho r\leq1$.
Consequently, by combining the inequalities \eqref{liu03-e} and \eqref{liu03-f}, we get
\begin{eqnarray}\label{HLP-eq3e}
B_H(f,r) &=& r^m\left (|a_{0}|+|b_{0}|+\sum_{k=1}^{\infty} |a_{k}|r^{pk}+\sum_{k=1}^{\infty }|b_{k}|r^{pk}\right )\nonumber\\
&\leq& 2r^m\left (a+ \frac{1+d}{2}\,\frac{r^{p}(1-a^{2})}{\sqrt{1-a^{2}r^{p}\rho^{p}}}\frac{1}{\sqrt{1-\rho^{-p}r^{p}}}\right ).
\end{eqnarray}
We wish to maximize the right hand side of above. For this, we need to consider the cases $a\geq r^{p}$ and $a<r^{p}$, separately.
Note that our choice of $\rho$ is such that $\rho r\leq 1$.
\bca Assume that $a\geq r^{p}$.
\eca
In this case we set $\rho=\frac{1}{\sqrt[p]{a}}$ and obtain from \eqref{HLP-eq3e} that
\be \label{liu04}
B_H(f,r)\leq 2r^{m}\psi(a) ~\mbox{ for $a\geq r^{p}$},
\ee
where we let $\alpha=r^{p}$ and
\begin{equation*}
\psi(x)=x+\frac{\alpha (1+d)}{2}\cdot \frac{1-x^{2}}{1-\alpha x},\quad  x\in[0,1].
\end{equation*}
Simple computation shows that, when $\alpha\geq\frac{1}{2+d}$, $\psi(x)$ attains its maximum at $x=x_{1}$, where
\begin{equation*}
x_{1}=\left (1-\sqrt{\frac{1+d}{3+d}}\sqrt{1-\alpha^{2}}\right )\frac{1}{\alpha},
\end{equation*}
and thus, $\psi(x)\leq\psi(x_{1})$. On the other hand, when $\alpha<\frac{1}{2+d}$, $\psi(x)$ is monotonically increasing
for $x\in[0,1]$ so that $\psi(x)\leq\psi(1)=1$.
Consequently, for the $r=r_{p,m,d}$ defined as in Theorem \ref{HLP-th2-new} and $\alpha\geq\frac{1}{2+d}$, it follows
from \eqref{liu04}  that
\begin{equation}
B_H(f,r) \leq 2r^m\psi(x_{1})=\frac{2}{r^{p-m}}\left (2+d-\sqrt{(1+d)(3+d)}\sqrt{1-r^{2p}}\right )  =1,
 \label{liu05}
\end{equation}
where we have used Lemma \ref{HLP-lem2} for the equality sign on the right. When $\alpha<\frac{1}{2+d}$, we have
 $\psi(x)\leq\psi(1)=1$ and thus,
\begin{equation*}
B_H(f,r)\leq 2r^{m}\psi(1)=2r^{m}.
\end{equation*}

\bca Assume that $a<r^{p}$.
\eca

In this case we set $\rho=\frac{1}{r}$ and obtain from \eqref{HLP-eq3e} that
\begin{equation}
B_H(f,r)\leq 2r^{m}\left (a+\frac{1+d}{2}\cdot r^{p}\frac{\sqrt{1-a^{2}}}{\sqrt{1-r^{2p}}}\right )
\leq (3+d)r^{p+m}\leq 1.
\label{liu06}
\end{equation}
Here the second inequality on the right follows from the argument that we omitted the critical point
$$a=\frac{\sqrt{1-r^{2p}}}{\sqrt{1+\left (\frac{(1+d)^2}{4}-1\right)r^{2p}}}
$$
because it is less than $r^{p}$ only in the case $r^{2p}>\frac{2}{3+d}\geq\frac{1}{2}$, which contradicts with Lemma \ref{HLP-lem1}.
The third inequality on the right in (\ref{liu06}) follows from Lemma \ref{HLP-lem1}.

Therefore, in both cases, for all $a\in [0,1)$, we have

\begin{enumerate}
\item[(i)]  when  $\alpha=r^{p}\geq\frac{1}{2+d}$,
\begin{equation*}
  B_H(f,r)\leq 1  ~\mbox{ for }~r\leq r_{p,m,d},
\end{equation*}
where $r_{p,m,d}$ is defined  as in Theorem \ref{HLP-th2-new}.

\item[(ii)]  when $\alpha=r^{p}<\frac{1}{2+d}$, since $\max\{2r^{m},(3+d)r^{p+m}\}=2r^{m}$, we have
\begin{equation*}
    B_H(f,r)\leq 2r^{m}.
\end{equation*}
\end{enumerate}
In summary, if $r^{m}=\frac{1}{2}$ and $r^{p}<\frac{1}{2+d}$, that is, if $\frac{p}{m}>\log_2 (2+d)$,
we have by the second case above
\begin{equation*}
 B_H(f,r)\leq 2r^{m}\leq 1,~\mbox{ for }~ r\leq\sqrt[m]{\frac{1}{2}}.
\end{equation*}
The extremal function for the case $d=1$ is $f(z)=h(z)+\overline{\lambda h(z)}$ with $h(z)=z^m$ and $|\lambda|=1$.

When $1\leq\frac{p}{m}\leq\log_2 (2+d)$, we apply the first case above to obtain
\begin{equation*}
 B(f,r) \leq \frac{1}{2} ~\mbox{ for }~r\leq r_{p,m,d},
\end{equation*}
where $r_{p,m,d}$  is defined  as in Theorem \ref{HLP-th2-new}.

For the case $d=1$, sharpness follows if we consider $f(z)=h(z)+\overline{\lambda h(z)}$ with $|\lambda|=1$,
$$h(z)=z^{m}\left (\frac{z^{p}-a}{1-az^{p}}\right ),\quad
a=\left (1-\frac{\sqrt{1-r_{p,m,1}^{2p}}}{\sqrt{2}}\right )\frac{1}{ r_{p,m,1}^{p}},
$$
and then calculate the Bohr radius for it. It coincides with $r$. \hfill $\Box$

\subsection{Proof of Theorem \ref{HLP-th1}}
Without loss of generality, we may assume that
$$
\max\{\|h\|_{\infty},\|g\|_{\infty}\}=1.
$$
\bca
Assume that $1\leq\frac{p}{m}\leq\log_{2}3$.
\eca
It follows from Corollary \ref{HLP-cor2}(1) and the hypothesis that $B(h,r) \leq \frac{1}{2}$ and $B(g,r) \leq \frac{1}{2}$ for $r\leq r_{p,m}$,
where $r_{p,m}$ is as in Theorem \ref{HLP-th1}. Adding these two inequalities shows  that
\begin{equation*}
B_H(f,r)=B(h,r) + B(g,r)\leq1 ~\mbox{ for $r\leq r_{p,m}$}.
\end{equation*}

\bca Assume that $\frac{p}{m}>\log_{2}3$.
\eca
 Applying the method of the previous case, Corollary \ref{HLP-cor2}(2) gives
\begin{equation*}
B_H(f,r)=B(h,r) + B(g,r)\leq 1,~~\mbox{for}~~r\leq\sqrt[m]{\frac{1}{2}}.
\end{equation*}

The extremal functions given in the statement are easy to verify.
\hfill $\Box$

\section{Bohr-type inequalities for harmonic mappings with a multiple zero at the origin}

\subsection{Proof of Theorem \ref{HLP-th3}}
By assumption, $|g'(z)|\leq d|h'(z)|$ for some $d\in [0,1)$. Then $\omega_{f}=\frac{g'}{h'}$ is analytic in the punctured
disk $0<|z|<1$ and has removable singularity at the origin with
$$\lim _{z\ra 0}\omega_{f}(z)=\frac{b_{pk+m}}{a_{pk+m}}
$$
so that $|b_{pk+m}|\leq d|a_{pk+m}|< |a_{pk+m}|$.

Since $h\in\mathcal{B}_{pk+m}$, by applying Lemma \ref{HLP-lem6} to the function $H(z)=\sum_{n=0}^{\infty }A_nz^{pn}$, $A_n =a_{p(n+k)+m}$,
one has
\be
\sum_{n=0}^{\infty}|A_n| r^{p n}+\left(\frac{1}{1+|A_0|}+\frac{r^{p}}{1-r^{p}}\right) \sum_{n=1}^{\infty}|
A_n |^{2} r^{2p n} \leq
|A_0|+\left(1-|A_0|^{2}\right) \frac{r^{p}}{1-r^{p}}.
\label{liu12}
\ee
Multiplying both sides of the inequality \eqref{liu12}  by the number $r^{p k+m}$, and then adding the term
$r^{p k+m}|A_0|^{2}\left(\frac{1}{1+|A_0|}+\frac{r^{p}}{1-r^{p}}\right)$ to both sides, we have
\begin{eqnarray} \hspace{1cm}
M_{pk+m}(h,r)&\leq&
{r^{p k+m}\left[|A_0|+ \frac{r^{p}}{1-r^{p}}+\frac{|A_0|^{2}}{1+|A_0|}\right]}  
= r^{p k+m}\left (\frac{r^{p}}{1-r^{p}}+ G(|A_0|) \right ),
\label{liu502}
\end{eqnarray}
where $G(t)= t +t^{2}(1+t)^{-1}.$ Since $G'(t)>0$ on $[0,1]$, it follows that $G(t)\leq G(1)=3/2$ and thus, \eqref{liu502} implies
\be\label{HLP-eq4}
 M_{pk+m}(h,r)\leq r^{p k+m}\left (\frac{r^{p}}{1-r^{p}}+\frac{3}{2}\right )=r^{p k+m}\left (\frac{3-r^{p}}{2(1-r^{p})}\right ).
\ee

Next, since $|g'(z)|\leq d|h'(z)|$ for $z\in \ID$,   we have (cf. \cite{AlKayPON-19})
\begin{eqnarray}
\sum_{n=k}^{\infty}\left|b_{p n+m}\right| r^{p n+m}\leq d \sum_{n=k}^{\infty}\left|a_{p n+m}\right| r^{p n+m}\quad\mbox{ for } r\leq 1/\sqrt[p]{3},
\label{liu503}
\end{eqnarray}
and, as in the proof of Theorem \ref{HLP-th2-new},
\begin{eqnarray}
\sum_{n=k}^{\infty}\left|b_{p n+m}\right|^2 r^{p(2 n-k)+m}\leq d^2\sum_{n=k}^{\infty}\left|a_{p n+m}\right|^2 r^{p(2 n-k)+m}.
\label{liu504}
\end{eqnarray}
Thus we conclude from (\ref{liu502}), (\ref{liu503}) and (\ref{liu504}) that
\begin{eqnarray}
N_{pk+m}(g,r)&\leq& d\sum_{n=k}^{\infty}\left|a_{p n+m}\right| r^{p n+m}+d^2\left(\frac{1}{1+\left|a_{p k+m}\right|}+\frac{r^{p}}{1-r^{p}}\right) \sum_{n=k}^{\infty}\left|a_{p n+m}\right|^{2} r^{p(2 n-k)+m}\nonumber
\end{eqnarray}
which by combining with (\ref{liu502}) gives

 \vspace{8pt}
 \noindent
 $M_{pk+m}(h,r)+ N_{pk+m}(g,r)$
\begin{eqnarray*}
 &\leq&
(1+d)\sum_{n=k}^{\infty}\left|a_{p n+m}\right| r^{p n+m}
+(1+d^2)\left(\frac{1}{1+\left|a_{p k+m}\right|}+\frac{r^{p}}{1-r^{p}}\right) \sum_{n=k}^{\infty}\left|a_{p n+m}\right|^{2} r^{p(2 n-k)+m}\\
&=& (d-d^2)\sum_{n=k}^{\infty}\left|a_{p n+m}\right| r^{p n+m} +(1+d^2)M_{pk+m}(h,r)\\
&\leq & (d-d^2)M_{pk+m}(h,r) +(1+d^2)M_{pk+m}(h,r) = (d+1)M_{pk+m}(h,r),
\label{liu506}
\end{eqnarray*}
which, by \eqref{HLP-eq4}, is less than or equal to $1$ if
$$r^{p k+m}\left (\frac{3-r^{p}}{2(1-r^{p})}\right ) \leq \frac{1}{1+d}, ~\mbox{ i.e., }~t_{k}(r) \geq 0,
$$
where $t_{k}(r)$ is given by \eqref{HLP-eq5}; that is,
\begin{equation*}
t_{k}(r)=\frac{2}{d+1}(1-r^{p})-r^{p k+m}\left(3-r^{p}\right).
\end{equation*}
This proves the first part of the assertion.

Now we prove the uniqueness of the solution in $(0, 1)$ of $t_{k}(r)=0$, we compute
that $t_{k}(0)= 2/(d+1)>0$, $t_{k}(1)=-2<0$, and
\begin{eqnarray*}
t_{k}'(r)&=&- \frac{2p}{d+1} r^{p-1}-3(p k+m) r^{p k+m-1}+(p(k+1)+m) r^{p(k+1)+m-1}\\
&=&-p\left( \frac{2}{d+1}r^{p-1}-r^{p(k+1)+m-1}\right)-(p k+m) r^{p k+m-1}\left(3-r^{p}\right)\\
&\leq &-p\left(r^{p-1}-r^{p(k+1)+m-1}\right)-(p k+m) r^{p k+m-1}\left(3-r^{p}\right)<0,
\end{eqnarray*}
showing that $t_{k}(r)$ is a decreasing function of $r$ in $(0, 1)$, and thus,
$t_{k}(r)=0$ has a unique root in $(0, 1)$.
\hfill $\Box$
\vskip 3mm

\subsection{Proof of Theorem \ref{HLP-th4}}
By assumption $h\in\mathcal{B}_{pk+m}$. Therefore, as in the proof of Theorem \ref{HLP-th3}, we can apply Lemma \ref{HLP-lem6} to the function $H(z)=\sum_{n=0}^{\infty }A_nz^{pn}$, $A_n =a_{p(n+k)+m}$. Thus, \eqref{liu12} holds.
 Multiplying the inequality  \eqref{liu12} by $r^{p k+m}$ gives
\begin{eqnarray*}
M_{pk+m}^{1}(h,r)&\leq & r^{p k+m}\left[|A_0|+\left(1-|A_0|^{2}\right) \frac{r^{p}}{1-r^{p}}\right].
\end{eqnarray*}

Now, we can maximize the right hand side with respect to $|A_0|$ by fixing $r$. A simple calculation shows that
we arrive at the maximum value $r^{p k+m}M(r)$ which is achieved at $|A_0|=1$, if $r\in\left[0, \frac{1}{\sqrt[p]{3}}\right]$, and
at $|A_0|=\frac{1-r^{p}}{2 r^{p}}$ in the remaining cases. Thus, we have the maximum value
$$\frac{r^{p(k-1)+m}\left(5 r^{2 p}-2 r^{p}+1\right)}{4\left(1-r^{p}\right)}
$$
and therefore,
\begin{eqnarray*}
M_{pk+m}^{1}(h,r)&\leq & r^{p k+m}\left[|A_0|+\left(1-|A_0|^{2}\right) \frac{r^{p}}{1-r^{p}}\right]\leq \frac{r^{p(k-1)+m}\left(5 r^{2 p}-2 r^{p}+1\right)}{4\left(1-r^{p}\right)}.
\end{eqnarray*}
Again, as $g\in\mathcal{B}_{pk+m}$, we have similarly the inequality
$$
M_{pk+m}^{1}(g,r)\leq r^{p k+m}\left[|B_0|+\left(1-|B_0|^{2}\right) \frac{r^{p}}{1-r^{p}}\right]
\leq \frac{r^{p(k-1)+m}\left(5 r^{2 p}-2 r^{p}+1\right)}{4\left(1-r^{p}\right)},
$$
where $B_0=b_{pk+m}$. Adding the two resulting inequalities yields that
\begin{eqnarray*}
M_{pk+m}^{1}(h,r) +M_{pk+m}^{1}(g,r) \leq \frac{r^{p(k-1)+m}\left(5 r^{2 p}-2 r^{p}+1\right)}{2\left(1-r^{p}\right)}.
\end{eqnarray*}
Hence the desired inequality (\ref{liu51}), i.e.
$M_{pk+m}^{1}(h,r) +M_{pk+m}^{1}(g,r)\leq 1$, holds whenever $L_{k}(r) \geq 0$, where
\begin{equation*}
L_{k}(r)=2\left(1-r^{p}\right)-r^{p(k-1)+m}\left(5 r^{2 p}-2 r^{p}+1\right).
\end{equation*}
This proves the first part of the assertion.

Next, we prove the uniqueness of the solution in $(0,1)$ of $L_{k}(r)=0$. In fact, note that $L_{k}(0)=2>0$, $L_{k}(1)=-4<0$, and
\begin{eqnarray*}
L_{k}^{\prime}(r)&=&-pr^{p-1}\left(2  -r^{p(k-2)+m}\right) -r^{p(k-1)+m-1}Q(r^p),
\end{eqnarray*}
where
$Q(x)=5(p(k+1)+m) x^{2} -2(p k+m)x+p k+m$. It follows that $Q(x)>0$,  because the discriminant of the function $Q$
is less than $0$. This gives that  $L_{k}^{\prime}(r)<0$ and hence, $L_{k}(r)=0$ has a unique root $\tau_{k}$ in $(0, 1)$.

Finally, we verify the sharpness of the upper bound $\tau_{k}$ for the Bohr radius. We consider the function
$f(z)=h(z) +\overline{\lambda h(z)}$, $|\lambda|=1$, where
\begin{equation*}
h(z)=z^{p k+m}\left (\frac{a-z^p}{1-a z^p}\right ),~a\in (0,1).
\end{equation*}
For this function, we obtain that
\begin{equation*}
M_{pk+m}^{1}(h,r)+ M_{pk+m}^{1}(g,r)=2r^{pk+m}\left [a+\frac{(1-a^2)r^p}{1-r^p}\right ]
\end{equation*}
which equals $1$ for $r=\tau_{k}$ and $a=\frac{1-\tau_{k}^p}{2 \tau_{k}^p}$. This completes the proof of Theorem  \ref{HLP-th4}
\hfill $\Box$
 \vskip 3mm

\subsection{Proof of Theorem \ref{HLP-th5}}
Let $h\in\mathcal{B}_{pk+m}$. Then \eqref{liu502} holds, that is,
\be
M_{pk+m}(h,r)
\leq  r^{p k+m}\left[ |A_0|+ \frac{r^{p}}{1-r^{p}} +\frac{|A_0|^2}{1+|A_0|}\right ]\leq  r^{p k+m}\left[\frac{3}{2}+\frac{r^{p}}{1-r^{p}}\right ],
\label{liu12-1}
\ee
where $A_0 =a_{pk+m}$. The second inequality holds because the function $T$ defined by
$$T(x)=x+ \frac{x^2}{1+x}
$$
is a monotonically increasing function of $x\in [0,1]$ so that $T(x)\leq T(1)=3/2$.

Similarly, with $B_0 =b_{pk+m}$, we have
\be
M_{pk+m}(g,r)
\leq   r^{p k+m}\left[ |B_0|+ \frac{r^{p}}{1-r^{p}} +\frac{|B_0|^2}{1+|B_0|}\right ]\leq  r^{p k+m}\left[\frac{3}{2}+\frac{r^{p}}{1-r^{p}}\right ].
\label{liu12-2}
\ee
where $|B_0|\in [0,1]$.

Combining \eqref{liu12-1} and \eqref{liu12-2}  leads to
\begin{eqnarray}
M_{pk+m}(h,r) +M_{pk+m}(g,r)
&\leq & r^{p k+m}\left (3+\frac{2 r^{p}}{1-r^{p}}\right ),
\label{liu54}
\end{eqnarray}
if and only if $w_{k}(r) \geq 0$, where
\begin{equation*}
w_{k}(r)=\left(1-r^{p}\right)-r^{p k+m}\left(3-r^{p}\right).
\end{equation*}
This proves the first part of the assertion of the theorem.

Next, to prove the uniqueness of the solution in $(0,1)$ of $w_{k}(r)=0$, it is sufficient to observe that
$w_{k}(0)=1>0$, $w_{k}(1)=-2<0$, and
\begin{eqnarray*}
w_{k}^{\prime}(r)=-pr^{p-1}\left(1 -r^{pk+m}\right)-(p k+m) r^{p k+m-1}\left(3-r^{p}\right)<0.
\end{eqnarray*}

Finally, it is easy to verify that the extremal function has the form $f(z)=h(z) +\overline{\lambda h(z)}$,
where $h(z)=z^{p k+m}$ and $|\lambda|=1$. This completes the proof of Theorem \ref{HLP-th5}.
\hfill $\Box$

\subsection{Proof of Theorem \ref{HLP-th6}}
The proof is essentially similar to the proof of Theorem \ref{HLP-th5}. At first, from (\ref{liu54}) and the assumption that
$|A_0|=|B_0|=a$, it is obvious that
the required inequality (\ref{liu55}) is true if
\begin{equation}
2r^{p k+m}\left(a+\frac{r^{p}}{1-r^{p}}+\frac{a^{2}}{1+a}\right)=2r^{p k+m}\left[ \frac{ a+2 a^{2}+r^{p}\left(1-2 a^{2}\right)}{(1+a)\left(1-r^{p}\right)}\right] \leq 1,
\label{liu56}
\end{equation}
which holds if and only if $V_{k, a}(r) \geq 0$, where
\begin{equation*}
V_{k, a}(r):=(1+a)\left(1-r^{p}\right)-2r^{p k+m}\left[2 a^{2}+a+r^{p}(1-2 a^{2})\right] .
\end{equation*}
This proves the first part of the assertion of the theorem.

Next, we can prove the uniqueness of the solution of $V_{k, a}(r)=0$ in $(0, 1)$. It is obvious that
$V_{k, a}(0)=1+a>0$ and $V_{k, a}(1)=-2(1+a)<0$. Furthermore,
\begin{eqnarray*}
V_{k, a}^{\prime}(r)=-p r^{p-1}\left[1+a+2\left(1-2 a^{2}\right) r^{p k+m}\right]-2 r^{p k+m-1}(p k+m)\left[2 a^{2}+a+\left(1-2 a^{2}\right) r^{p}\right]
\end{eqnarray*}
and it is easy to obtain that $V_{k, a}^{\prime}(r)<0$ and thus, $V_{k, a}(r)=0$ has the unique root $\varsigma_{k}=\varsigma_{k}(a)$ in
the interval $(0, 1)$.

Finally, we verify the sharpness of the upper bound $\varsigma_{k}$ for the Bohr radius. Consider
$f(z)=h(z) +\overline{\lambda h(z)}$, where
$$h(z)=z^{p k+m}\left (\frac{a-z^p}{1-a z^p}\right ), ~a\in [0,1],
$$
$|\lambda|=1$ and $a\in [0, 1]$ is fixed.
In this case, we get for the left hand side of (\ref{liu55}) (for simplicity call it as $W(r)$ ) takes the form
\begin{eqnarray*}
W(r)&=&2r^{pk+m}\left [a+\frac{(1-a^2)r^p}{1-r^p}\right]+2\left(\frac{1}{1+a}+\frac{r^{p}}{1-r^{p}}\right)a^{2} r^{p k+m}\\
&=&2r^{p k+m}\left[ \frac{ a+2 a^{2}+r^{p}\left(1-2 a^{2}\right)}{(1+a)\left(1-r^{p}\right)}\right].
\end{eqnarray*}
Comparison of this expression with the right hand side of the equation in formula (\ref{liu56}) delivers the asserted sharpness.
The proof of Theorem \ref{HLP-th6} is complete.\hfill $\Box$

\subsection*{Acknowledgments}
This research of the first two authors are partly supported by Guangdong Natural Science Foundations (Grant No. 2021A030313326).
The work of the third author was supported by Mathematical Research Impact Centric Support (MATRICS) of
the Department of Science and Technology (DST), India  (MTR/2017/000367).

%

\end{document}